\documentclass{article}

\usepackage{a4}

\usepackage[latin1]{inputenc}

\usepackage[T1]{fontenc}

\usepackage[english]{babel}

\usepackage{amsmath}

\usepackage{amssymb}

\usepackage{vmargin}

\setmarginsrb{2cm}{2.5cm}{2.5cm}{2.5cm}{0cm}{0.8cm}{0cm}{0.8cm}

\newcounter{proposition}[section]

\newcounter{definition}[section]

\newcounter{theoreme}[section]

\newcounter{lemme}[section]

\newcounter{corollaire}[section]

\makeatletter

\def\theproposition{\thesection.\@arabic\c@proposition}

\def\thedefinition{\thesection.\@arabic\c@definition}

\def\thetheoreme{\thesection.\@arabic\c@theoreme}

\def\thelemme{\thesection.\@arabic\c@lemme}

\makeatother

\newenvironment{enonce}[1]{\noindent{\textbf{\scshape{}#1}}---\,\,\begin{itshape}}{\end{itshape}}

\newenvironment{proposition}[1][]{\refstepcounter{proposition}

\begin{enonce}{Proposition \theproposition{}#1 }}{\end{enonce}}

\newenvironment{definition}[1][]{\refstepcounter{definition}

\begin{enonce}{Definition \thedefinition{}#1 }}{\end{enonce}}

\newenvironment{theoreme}[1][]{\refstepcounter{theoreme}

\begin{enonce}{Theorem \thetheoreme{}#1 }}{\end{enonce}}

\newenvironment{lemme}[1][]{\refstepcounter{lemme}

\begin{enonce}{Lemma \thelemme{}#1 }}{\end{enonce}}

\newenvironment{demo}[1][]{\noindent{\scshape{}Proof{}#1 }---\,\,}{\,\textbf{Q.E.D.}}

\newenvironment{corollaire}[1][]{\refstepcounter{corollaire}

\begin{enonce}{Corollary \thecorollaire{}#1 }}{\end{enonce}}

\renewcommand{\leq}{\leqslant}

\renewcommand{\geq}{\geqslant}

\DeclareMathOperator{\Card}{Card}

\newcommand{\N}{\mathbb{N}}

\newcommand{\Z}{\mathbb{Z}}

\title{The direct extension theorem}

\author{Joseph Ayoub}

\begin{document}

\maketitle

\begin{abstract}
\emph{The problem of group extension can be divided into two
  sub-problems. The first is to find all the possible extensions of
  $H$ by $K$. The second is to find the different ways a group $G$ can
  arise as an extension of $H$ by $K$. Here we prove that the direct
  product $H\times K$ can arise as an extension of $H$ by $K$ in an
  essentially unique way: that is the direct extension. 
I would like to thank Yacine Dolivet for drawing my attention to the
direct ``extension theorem'', and Anne-Marie Aubert as well as
Charles-Antoine Louet for their support.}
\end{abstract}

\medskip

\begin{tableofcontents}

\end{tableofcontents}

\bigskip

NOTATIONS:

\begin{itemize}
\item $1$:       \emph{is the identity group}
\item $G$, $H$, $K$, $L$... \emph{will denote finite groups} 
\item $H\leq G$: \emph{``$H$ is a subgroup of $G$''}
\item $H<G$: \emph{``$H$ is a subgroup of $G$ distinct from $G$''}
\item $H\unlhd G$: \emph{``$H$ is a normal subgroup of $G$''}
\item $\mathcal{Z}(G)$: \emph{is the center of $G$}
\item $G'$: \emph{is the derived group of $G$}
\item $\Card(G)$ \emph{or} $|G|$: \emph{is the order of the finite
    group $G$}
\end{itemize}

\pagebreak

\section{Statement of the theorem}
In this paper, we prove the following theorem, which we call the
``direct extension'' theorem.

\medskip

\begin{theoreme}
Let $G$, $H$ and $K$ be three finite groups. If $G$ and
$H\times K$ are isomoprhic, then every group extension  $1
\longmapsto H \longmapsto G \longmapsto K \longmapsto 1$ is a direct extension.

\end{theoreme}

\medskip

We may reformulate the theorem as follows.

\medskip

\begin{theoreme}
Let $G=H.K$ be a 
decomposition of the finite group $G$ into direct factors
and $H_0$ be a normal subgroup of $G$. Assume that $H_0$ and  $H$ are isomorphic, as
well as $G/H_0$ and $K$. Then $H_0$ is a direct factor of $G$ (that
is, there exists $K_0 \unlhd G$ such that $G=H_0.K_0$ and $H_0\cap K_0=1$).

\end{theoreme}

\medskip

We use the latter statement of the theorem in the proof. Assume that
the theorem does not hold, then there is counter-example
$(G=H.K;H_0)$, which is minimal with respect to $\Card G$. We shall
derive a contradiction from the existence of $G$.

\section{A few preliminary general results}

In this section, $G$ is any finite group, not necesseraly the group that
appears in the theorem.

\subsection{Subgroups of a direct product $G=H.K$}

We give here some useful simple results:

\medskip

\begin{proposition} 
Let $L$ be a subgroup of $G$, which we do not assume to be normal,
and $G=H.K$ a decomposition of $G$ into direct factors.
Assume that $H\subseteq L$. Then $L=H.(L\cap K)$. In particular, $H$
is a direct factor of $L$.
\end{proposition}

\medskip

\begin{demo}
It is clear that $H.(L\cap K) \subseteq L$. Now, every $l\in L$ may be
written $l=h.k$, with $h\in H$ and $k\in K$. $H\subseteq
L$, which shows that $h\in L$. Thus $k\in L\cap K$ and $l\in H.(L\cap K)$.
\end{demo}

\medskip

\begin{proposition}
Let $G=H.K$ be a decomposition of $G$ into direct factors. Then
$G'=H'.K'$ and $\mathcal{Z}(G)=\mathcal{Z}(H).\mathcal{Z}(K)$.

\end{proposition}

\medskip

\begin{demo}
The first assertion comes from the following formula
$[h_1.k_1,h_2.k_2]=[h_1,h_2].[k_1,k_2]$ which is true for $h_i\in H$
and $k_i\in K$.
The second assertion is trivial. 
\end{demo}

\subsection{Coprime direct factors of a finite group $G$}

In this paragraph, we make use of the famous Remak-Krull-Schmidt theorem on the 
decomposition of finite groups into indecomposable direct factors.

\medskip

\begin{definition}
Let $A$ and $B$ be two finite groups. $A$ and $B$ are said to be
coprime if no non-trivial direct factor of $A$ is isomorphic to a direct factor of
$B$.

\end{definition}

\medskip

\begin{proposition}
Let $A$ and $B$ be two direct factors of the finite group $G$. Assume
that $A$ and $B$ are coprime. Then $A\cap B=1$ and $A.B$ is
a direct factor of $G$.
\end{proposition}

\medskip

\begin{demo}
Let $C$ be a direct factor of $G$ such that:

\begin{itemize}
\item $A$ is a direct factor of $C$,
\item The indecomposable direct factors of $C$ are isomorphic to direct
factors of $A$,
\item $C$ is maximal having these properties.
\end{itemize}

To find such a $C$, it suffices to take a maximal element
of the set of the direct factors of $G$ satisfying the first two
conditions. That set is non empty because it contains $A$.

\pagebreak

Define also $D$, a direct factor of $G$, such that:

\begin{itemize}
\item $B$ is a direct factor of $D$,
\item The indecomposable direct factors of $D$ are not isomorphic to any
of the indecomposable factors of $A$ (this implies the same thing with $C$
instead of $A$).
\item $D$ is maximal having these properties.
\end{itemize}

As before, to build such a $D$, it suffices to take a maximal element
of the set of the direct factors of $G$ satisfying the two first
conditions. 
That set contains $B$, because $B$ and $A$ are coprime, as well as $B$
and $C$.

It is clear that $C$ and $D$ are coprime. Let us now prove that 
$G=C.D$, and that it is a direct product (i.e. $C\cap D=1$). We will
have finished proving the proposition because $A\cap B \subseteq C\cap
D=1$, and if
$C=A.U$, $D=B.V$ are
decompositions into direct factors
$G=(A.B).(U.V)$ is a 
decomposition
of $G$ into direct factors.

Let $G=A_1...A_m.B_1...B_n$, be a 
decomposition of $G$ into indecomposable direct factors
, such that each of the $A_i$ be isomorphic to a direct factor of
$A$, and such that none of the $B_j$ be isomorphic to a direct factor of $A$.
We now show that $C$ is isomorphic to $A_1...A_m$. We write
$G=C.S$ as a direct product. It is clear that every indecomposable
direct factor of $S$ is, up to isomorphism, one of the $B_j$, by
maximality of $C$. According to the Remak-Krull-Schmidt theorem, $S$
is isomorphic to a partial product of $A_1...B_n$, therefore it is
isomorphic to a direct factor of $B_1...B_n$. Thus, $|S|\leq
|B_1...B_n|$ and $|C|\geq
|A_1...A_m|$. Again using Remak-Krull-Schmidt, $C$ is isomorphic to a
direct factor of $A_1...A_m$. It follows that $C$ is isomorphic to
$A_1...A_m$.

Likewise, we show that $D$ is isomorphic to $B_1...B_n$. In particular
$|C|\times |D|=|G|$. We still have to prove that $C\cap D=1$. Consider
the subgroup $C.D$. It contains $C$ as a direct factor.
But $C$ is isomorphic to $A_1...A_m$. So $C.D \simeq
A_1...A_m\times((C.D)/C)$. Now $D$ is also a direct factor of
$C.D$, so according to the Remak-Krull-Schmidt theorem, $D\simeq B_1...B_n$ is
isomorphic to a direct factor of $(C.D)/C\simeq D/(C\cap D)$. Thus $|D|\leq
|D/(C\cap D)|$ which shows that  $C\cap D=1$.
\end{demo}

\medskip

\begin{corollaire}
Let $G=B.C$ be a decomposition of $G$ into direct factors. 
Let $A$ be a direct factor of $G$, such that $A$ and $B$ are coprime.
Then the projection of $A$ onto $C$ is a direct factor of $G$.
\end{corollaire}

\medskip

\begin{demo}
We know that $B.A$ is a direct factor of $G$. But $B.A=B.p_C(A)$, where
$p_C$ is the projection onto $C$ with respect to $B$. This shows that
$p_C(A)$ is a direct factor of $B.A$, hence also of $G$.
\end{demo}

\subsection{Directly decomposable subgroups of $G$}

In this paragraph, we define the concept of directly decomposable
subgroups, and show two propositions that will be needed later on.

\medskip

\begin{definition}
Let $G$ be a finite group and $D$ a subgroup of $G$. $D$ is said
to be direcly decomposable in $G$ if $D=(H\cap D).(K\cap D)$
for every decomposition $G = H.K$ of $G$ into direct factors.
\end{definition}

\medskip

\begin{proposition}
If $D$ is directly decomposable in $G$, then for every decomposition
$G=H_1...H_m$ of $G$ into direct factors we have
$D=(H_1\cap D)...(H_m\cap D)$. Furthermore, if $D$ is a normal subgroup of
$G$, $G/D=((H_1.D)/D)...((H_m.D)/D)$ is a decomposition of $G/D$ into
direct factors.

\end{proposition}

\medskip

\begin{demo}
We start with the proof of the first part of the statement. Consider a fixed
decomposition $G=H_1...H_m$ of $G$ into direct factors .
It is easy to see that $D=(H_1\cap D)...(H_m\cap D)$ is equivalent to the
following statement : for all $d=h_1...h_m \in D$ with $h_i\in H_i$, $h_i\in
D$.
Let $d=h_1...h_m\in D$. As $D$ is directly decomposable
in $G$, $D=(H_i\cap D).((H_1...H_{i-1}.H_{i+1}...H_m)\cap D)$
for all $1\leq i\leq m$. But then $h_i\in H_i\cap D$.

Assume now that $D$ is normal in $G$. We have
$G/D=((H_1.D)/D)...((H_m.D)/D)$. Since $(H_i.D)/D\simeq H_i/(H_i\cap
D)$, it is true that $$|(H_1.D)/D|\times ...\times |(H_m.D)/D|=(\frac{|H_1|}{|H_1\cap
D|})\times ...\times (\frac{|H_m|}{|H_m\cap D|})=|G|/|D|=|G/D|,$$

which shows that $G/D=((H_1.D)/D)...((H_m.D)/D)$ is necessarily a decomposition
of $G/D$ into direct factors. This completes the proof of our statement.
\end{demo}

\medskip

\begin{proposition}
If $T$ is a normal subgroup of $G$ such that
$T'=T\cap G'$, then $T'$ is directly decomposable in $G$.

\end{proposition}

\medskip

\begin{demo}
Let $G=L.M$ be a decomposition of $G$ into direct factors. Write
$T=(\overline{A},A;\overline{B},B)_{\varphi}$. Clearly $T\subseteq
\overline{A}.\overline{B}$, so $T'\subseteq
\overline{A}'.\overline{B}'$. But $T$ is a normal subgroup of $G$ so
according to paragraph 2.1 $\overline{A}'\subseteq A$ and
$\overline{B}'\subseteq B$. Thus, using
$A.B\subseteq T$, we get the following chain of inclusions:
$$T'\subseteq \overline{A}'.\overline{B}'\subseteq (A.B)\cap G
\subseteq T\cap G'=T'.$$

These inclusions are then necessarily equalities, so we have
$T'=\overline{A}'.\overline{B}'=(T'\cap L).(T'\cap M)$.
This shows that the subgroup is directly decomposable. 
\end{demo}

\section{Two special cases of the theorem}

In this section, we prove the theorem in the two special cases:
\begin{itemize}
\item $G$ is a commutative group,
\item $G'$, the derived subgroup of $G$, is equal to $G$.
\end{itemize}

As in part 2, $G$ stands for any finite group.

\subsection{The case of commutative groups}

We show the following result:

\medskip

\begin{proposition}
Let $G$ be a commutative finite group and $G=H.K$ a decomposition of $G$
into direct factors. Let $H_0$ be a subgroup of $G$, such that $H$ and
$H_0$ are isomorphic, as well as $G/H_0$ and $K$. Then $H_0$ is
a direct factor of $G$.

\end{proposition}

\medskip

As in the general case, we prove the statement \emph{ab absurdo} taking
a counter-example $(G=H.K,H_0)$ which is minimal with respect to $\Card G$.

We shall need a few lemmas, which we now state:

\medskip

\begin{lemme}
$H_0$ does not contain a non-trivial direct factor of $G$.
\end{lemme}

\medskip

\begin{demo}
Let $L\subseteq H_0$ be a non-trivial direct factor of $G$.
$L$ is also a direct factor of $H_0$. After changing
$H$ if needed, we may assume that $L\subseteq H$ because, since
$H\simeq H_0$, $H$ also contains a direct factor isomorphic to $L$.

Let us now consider the group $G/L$, which has striclty less elements
than $G$. $G/L \simeq H/L \times K$, $H/L \simeq H_0/L$ and
$(G/L)/(H_0/L) \simeq G/H_0 \simeq K$. Using minimaliy of $G$,
we may take $P$ a subgroup of $G$, $L \subset P$, such that
$G/L=(H_0/L).(P/L)$ is a decompostion into direct factors.
We then have $G=H_0.P$, with $H_0\cap
P=L$.

But $L$ is a direct factor of $P$, so $P=L.K_0$. This yields
$G=H_0.K_0$, with $H_0\cap K_0=1$. But that is impossible, because of our
assumptions on $G$.
\end{demo}

\medskip

\begin{lemme}
$G$ is a commutative $p$-group.
\end{lemme}

\medskip

\begin{demo}
Write $G=S_p.L$ where $S_p$ is a non-trivial $p$-Sylow of $G$
and $L$ is the Hall $p'$-group of $G$. We would like to show that $L=1$.
Clearly, $S_p=(S_p\cap H).(S_p\cap
K)$ and $S_p\cap H \simeq S_p\cap H_0$. 
Moreover $S_p/(S_p\cap H_0)$ and
$(S_p.H_0)/H_0$ are isomorphic.
But the latter group is isomorphic to the $p$-Sylow of $G/H_0$, which is
therefore isomorphic to $S_p\cap K$.

Now if $L>1$, minimality of $G$ shows that $S_p\cap H_0$
is a direct factor of $S_p$, therefore also of $G$. This contradicts
lemma 3.1, so we do have $G=S_p$.
\end{demo}

\medskip

\begin{definition}
Let $A$ be a group. We write $\omega(A)$ for the  l.c.m. of
the orders of  the elements of $A$. 
If $n$ is an integer, let $\mho^n(A)$ be the smallest subgroup of
$A$ which contains $x^n$ for all $x\in A$.
\end{definition}

\pagebreak

\begin{lemme}
$\omega(H) \geq \omega(K)$.

\end{lemme}

\medskip

\begin{demo}
Again, we proceed \emph{ab absurdo}. Assume $\omega (K) > \omega(H)$.
We know that $G$ is a commutative $p$-group, so $\omega (K)=p^k$
and $\omega (H)=p^h$, with $h\leq k-1$.
But then $\mho^{p^{k-1}}(H)=\mho^{p^{k-1}}(H_0)=1$ and
$\mho^{p^{k-1}}(G)=\mho^{p^{k-1}}(K)>1$.
We have $\mho^{p^{k-1}}(G/H_0)\simeq \mho^{p^{k-1}}(K)$, but
$\mho^{p^{k-1}}(G/H_0)=(\mho^{p^{k-1}}(G).H_0)/H_0 \simeq
\mho^{p^{k-1}}(G)/(H_0\cap \mho^{p^{k-1}}(G))$. It follows that
$H_0\cap \mho^{p^{k-1}}(G)=1$.

On the other hand,
$G/\mho^{p^{k-1}}(G)=((H.\mho^{p^{k-1}}(G))/\mho^{p^{k-1}}(G)).((K.\mho^{p^{k-1}}(G))/\mho^{p^{k-1}}(G))$
is a decomposition into direct factors,
$(H_0.\mho^{p^{k-1}}(G))/\mho^{p^{k-1}}(G)\simeq H_0$ and
$(H.\mho^{p^{k-1}}(G))/\mho^{p^{k-1}}(G)\simeq H$.

Finally, 
$(G/\mho^{p^{k-1}}(G))/((H_0\mho^{p^{k-1}}(G))/\mho^{p^{k-1}}(G))\simeq
(G/H_0)/((H_0\mho^{p^{k-1}}(G))/H_0)\simeq K/\mho^{p^{k-1}}(K)$.
Therefore, again by minimality of $G$, we may consider
a subgroup $K_0$ of $G$ which contains $\mho^{p^{k-1}}(G)$ and such that
$G/\mho^{p^{k-1}}(G)=(H_0.\mho^{p^{k-1}}(G))/(\mho^{p^{k-1}}(G)).(K_0/\mho^{p^{k-1}}(G))$
is a decomposition into direct factors. But then $G=H_0.K_0$ and
$\Card K_0=\Card K$. So $H_0$ is a direct factor of $G$. We have reached a contradiction.
\end{demo}

\medskip

\begin{lemme}
Let $A$ be a commutative $p$-group. Let $D$ be a cyclic subgroup of
$A$ with maximal order. $D$ is then a direct factor of $A$.
\end{lemme}

\medskip

\begin{demo}
We prove the lemma by induction on the order of $A$. Write
$A=B.C$, where $B$ is an indecomposable direct factor and $C$  a
supplementary normal subgroup to $B$ in $A$. As $D$ is cyclic, either $B\cap D=1$ or
$C\cap D=1$. If $C\cap D=1$ then as $\Card D \geq \Card B$, $A=C.D$ is a
decomposition into direct factors.

From now on, we assume $B\cap
D=1$. $q$, the projection onto $C$ according to $B$ is
a monomorphism on $D$. Thus $q(D)$ is a cyclic subgroup of maximal
order of $C$. The induction hypotheses then shows that $q(D)$ is a 
direct factor of $C$. So $B.D=B.q(D)$ is a direct factor
of $A$, therefore $D$ is a direct factor of $A$.
\end{demo}

\medskip

We now prove the actual proposition.

\begin{demo}
Let $D$ be an indecomposable direct factor of $H_0$, of maximal order.
According to lemma 3.3, $D$ is a cyclic subgroup of $G$ of maximal order.
According to lemma 3.4, $D$ is also a direct factor of $G$. But
that contradicts lemma 3.1. The proposition is true.
\end{demo}

\medskip

\subsection{The case where $G'=G$}

\begin{proposition}
Let $G$ be a finite group such that $G$ is equal to the derived group $G'$
and $G=H.K$ a decomposition of $G$ into direct factors.
Let $H_0 \unlhd G$
be isomorphic to $H$ and such that $G/H_0$ is isomorphic to $K$. Then $H_0$
is a direct factor of $G$.
\end{proposition}

\medskip

\begin{demo}
As usual, consider a minimal counter-example $(G=H.K;H_0)$.

We show that $H_0$ is directly decomposable in $G$.
$H_0'=H_0$ because $H_0$ is isomorphic to $H$ and
$G'=H'.K'=H.K$. It follows that $H_0'=H_0\cap G=H_0\cap G'$.
Proposition 2.5. then gives the result.

Moreover, $H_0$ does not contain a non-trivial direct factor of $G$.
It is an exercise to show that, proceeding in rather the same way
as in lemma 3.1.

Now, let $G=H_1...H_m.K_1...K_n$ be a decomposition of $G$ into indecomposable
direct factors such that $H=H_1...H_m$ and $K=K_1...K_n$. As
$H_0$ is directly decomposable in $G$, proposition 2.4 shows that
$$
G/H_0\simeq
(H_1/(H_1\cap H_0))\times ...\times (H_m/(H_m\cap H_0)) \times (K_1/(K_1\cap
H_0))\times ...\times (K_n/(K_n\cap H_0)).
$$

But none of the $H_i/(H_i\cap H_0)$ and none of the $K_j/(K_j\cap H_0)$ are
trivial. So $G/H_0 \simeq K$ contains at least $n+m$
indecomposable direct factors in a decomposition into irreducible direct
factors. We deduce that $n+m\leq n$, and
$m\leq 0$. We have reached a contradiction and out proposition is proved.
\end{demo}

\medskip

We have shown that if $(G=H.K;H_0)$ is a counter-example to the theorem, then
$1<G'<G$. It is the starting point of our proof of the theorem.

\pagebreak

\section{A few preliminary lemmas}

From now on, $(G=H.K;H_0)$ is our minimal counter-example to
the theorem. A few lemmas follow, which are useful to describe $H_0'$ and
$\mathcal{Z}(H_0)$ in $G$.

\medskip

\begin{lemme}
We have the following properties.

\begin{enumerate}
\item $H_0'=H_0\cap G'$,
\item There exists $M\unlhd G$ such that $G=M.H_0$ and $M\cap H_0=H_0'$,
\item $G/H_0'$ and $G/H'$ are isomorphic, as well as $M/H_0'$ and $K$,
\item $G/H_0'=(H_0/H_0').(M/H_0')$ is a decomposition of $G/H_0'$ into direct
factors.
\end{enumerate}

\end{lemme}

\medskip

In what follows, we fix a subgroup $M$ once and for all, which complies
with point 3 of lemma 4.1.We now prove the lemma.

\medskip

\begin{demo}
Let us start with point 1. Clearly $G'=H'.K'$ is a decomposition of $G'$ into
direct factors. Therefore $G'/H' \simeq K'\simeq (G/H_0)'$. But
$(G/H_0)'=(G'.H_0)/H_0\simeq G'/(H_0\cap G')$. It follows that
$\Card H_0\cap G'=\Card H'= \Card H_0'$. But $H_0'\subseteq H_0\cap
G'$, so the two subgroups are actually equal.

To show point 2, notice two things. First, $G/G'=(HG'/G').(KG'/G')\simeq
(H/H')\times (K/K')$. Then $(G/G')/((H_0.G')/G')\simeq
(G/(H_0.G'))\simeq ((G/H_0)/(G'.H_0)/H_0)= (G/H_0)/(G/H_0)'\simeq
K/K'$. As $G/G'$ is commutative, we may use part 3
to show that there exists a normal subgroup $M$ of $G$
containing $G'$ such that $G/G'=((H_0.G')/G').(M/G')$ is
a decomposition of $G/G'$ into direct factors. This proves point 2.

We now proceed to prove the two last statements.
$G/H_0'=(H_0/H_0').(M/H_0')$ is a decomposition into direct factors.
Moreover $M/H_0'\simeq (G/H_0')/(H_0/H_0')\simeq
G/H_0\simeq K$, and $H_0/H_0'\simeq H/H'$. So $G/H_0'$ and $G/H'$ are isomorphic.
\end{demo}

\medskip

We now state a corollary of the above lemma, which is crucial in the
proof of the "direct extension" theorem.

\medskip

\begin{corollaire}
$H_0'$ is directly decomposable in $G$.

\end{corollaire}

\medskip

\begin{demo}
It is an immediate consequence of proposition 2.5 and the above lemma.
\end{demo}

\medskip

The corollary shows that taking the quotient by $H_0'$ is compatible
with any decomposition of $G$ into direct factors. More precisely,
if $G=L.M$ is a decomposition of $G$ into direct factors then
$G/H_0'=((L.H_0')/H_0').((M.H_0')/H_0')$ is also a decomposition
into direct factors.

\medskip

\begin{lemme}
We have the following properties:
\begin{enumerate}
\item $\mathcal{Z}(H_0)=H_0\cap \mathcal{Z}(G)$,
\item $\mathcal{Z}(H_0)$ is a direct factor of $\mathcal{Z}(G)$,
\item $\mathcal{Z}(M/H_0')=(M\cap (\mathcal{Z}(G).H_0))/H_0'$.
\end{enumerate}

\end{lemme}

\medskip

\begin{demo}
$\mathcal{Z}(H_0)\simeq \mathcal{Z}(H)$ and $\mathcal{Z}(K)\simeq
\mathcal{Z}(G)/\mathcal{Z}(H)$. But $\mathcal{Z}(K)\simeq
\mathcal{Z}(G/H_0)\supseteq (\mathcal{Z}(G).H_0)/H_0\simeq
\mathcal{Z}(G)/(H_0 \cap \mathcal{Z}(G))$. Therefore $|H_0\cap
\mathcal{Z}(G)|\geq |\mathcal{Z}(H)|=|\mathcal{Z}(H_0)|$. But we know
that $H_0\cap \mathcal{Z}(G) \subseteq \mathcal{Z}(H_0)$. Hence the
equality of the two groups. We have also achieved
$\mathcal{Z}(G/H_0)=(\mathcal{Z}(G).H_0)/H_0$.
This completes the proof of the first point.

For point number 2, notice that $\mathcal{Z}(G)/\mathcal{Z}(H_0)=
\mathcal{Z}(G)/(H_0\cap \mathcal{Z}(G))\simeq
(\mathcal{Z}(G).H_0)/H_0=\mathcal{Z}(G/H_0)$. As $G$ is not equal to 
its center, it follows from the minimality of $G$ that
$\mathcal{Z}(H_0)$ is a direct factor of $\mathcal{Z}(G)$.

To prove the 3rd point, consider the natural isomorphism $\sigma
:M/H_0' \longmapsto G/H_0$. We have
$$ \mathcal{Z}(M/H_0')=\sigma^{-1}((\mathcal{Z}(G).H_0)/H_0)=\{x.H_0'\in
M/H_0'/ x.H_0 \in (\mathcal{Z}(G).H_0)/H_0\} $$

$$ =\{x.H_0'/x\in
\mathcal{Z}(G).H_0,  x\in M\}= ((\mathcal{Z}(G).H_0)\cap M)/H_0'. $$

Hence the announced result.
\end{demo}

\medskip

\begin{lemme}
$H_0$ does not contain a direct factor of $G$ other than 1. Similarly,
$H_0$ is not contained in a direct factor of $G$ other than
$G$.

\end{lemme}

\medskip

\begin{demo}
The first statement is left as an exercise. We prove the second one, which
is as simple as the first one.

Let $G=L.N$ be a decomposition of $G$ into direct factors. Assume
$N>1$ and $H_0\subseteq L$. Clearly $G/H_0=(L/H_0).((N.H_0)/H_0)$ is
a decomposition of $G/H_0$ into direct factors. Since $K\simeq
G/H_0$, $N\simeq
(N.H_0)/H_0$ is isomorphic to a direct factor of $K$. We may therefore
assume that $N\subset K$. We may then write $G/N\simeq H\times K/N$. But
$H\simeq (H_0.N)/N$, and $(G/N)/((H_0.N)/N)\simeq
(G/H_0)/((N.H_0)/H_0)\simeq K/N$, because $(N.H_0)/H_0$ is a direct factor
of $G/H_0$ which is isomorphic to $N$. It follows that $(H_0.N)/N$ is a
direct factor of $G/N$, using minimality of $G$.
Thus we have a normal subgroup $P$ of $G$ containing
$N$ such that $G=(H_0.N).P$ with $(H_0.N)\cap P=N$. It is now clear that
$G=H_0.P$ is a decomposition of $G$ into direct factors.
That contradicts our assumption on $G$.
\end{demo}

\section{The proof of the theorem}

We may now proceed with the actual proof of our theorem.

\medskip

\begin{proposition}
$H$ contains no non-trivial commutative direct factor.

\end{proposition}

\medskip

\begin{demo}
Again \emph{ab absurdo}. Let $A$ be a non-trivial commutative direct factor
of $H_0$, which exists since $H\simeq H_0$. We know that
$G/H_0'=(H_0/H_0').(M/H_0')$ is a decomposition of $G/H_0'$ into
direct factors. It is easy to see that $(A.H_0')/H_0'$ is a direct factor
of $H_0/H_0'$ and hence also of $G/H_0'$. Thus there exists
a direct factor $N/H_0'$ of $G/H_0'$ which contains $M/H_0'$ and is 
a supplementary subgroup of $(A.H_0')/H_0'$. But $G=A.N$ is then a 
decomposition of $G$ into direct factors since $(A.H_0')\cap N=H_0'$ and
it follows that $A\cap N \subseteq A\cap H_0'=1$ (since $A$ is a commutative
direct factor of $H_0$). We have shown that $A\subset H_0$ is also a direct
factor of $G$. This contradicts lemma 4.3.
\end{demo}

\medskip

\begin{proposition}
If $L$ is a non-commutative direct factor of $G$, then $L\cap H_0'>1$.

\end{proposition}

\medskip

\begin{demo}
\emph{Ab absurdo}. Let $L$ be a non-commutative direct factor of
$G$ such that $L\cap H_0'=1$. We may suppose that $L$
is indecomposable.  Then as $H_0'$ is directly decomposable in $G$,
$(L.H_0')/H_0'$ is a direct factor of $G/H_0'$ isomorphic to $L$. But
$L$ and $H_0/H_0'$ are coprime because $L$ is a non-commutative
indecomposable group and $H_0/H_0'$ is commutative. Then
$((L.H_0)'/H_0').(H_0/H_0')$ is a direct factor of $G/H_0'$ by
proposition 2.3.

Now let $H_0'\leq P\unlhd G$, such that
$G/H_0'=((L.H_0')/H_0').(H_0/H_0').(P/H_0')$ be a decomposition of
$G/H_0'$ into direct factors. Then we have $G=L.(P.H_0)$ and $L\cap
(P.H_0)= L\cap (L.H_0')\cap H_0=L\cap H_0'=1$. We have reached a
contradiction since this implies that $P.H_0$ is a direct factor of
$G$ distinct from $G$ and containing $H_0$.
\end{demo}

\medskip

\begin{proposition}
If $A$ is a commutative direct factor of $G$ then $H_0\cap A=1$
\end{proposition}

\medskip

\begin{demo}
Consider $H_0.A \subset G$. Clearly $\mathcal{Z}(H_0).A\subset \mathcal{Z}(G)$
and $\mathcal{Z}(H_0)$ is a direct factor of $\mathcal{Z}(G)$ (by lemma
4.2). Proposition 2.2 shows that there exists $B$ a supplementary
of $\mathcal{Z}(H_0)$ in $\mathcal{Z}(H_0).A$. It follows that
$H_0.A$=$H_0.B$ with $H_0\cap B=\mathcal{Z}(H_0)\cap B=1$. Thus $H_0$
is a direct factor of $H_0.A$. In the same way, $A$ is a direct factor
of $H_0.A$ because it is a direct factor of $G$. But according to 
proposition 5.1 $H_0$ and $A$ are coprime. Therefore $H_0\cap
A=1$, using once again proposition 2.3.
\end{demo}

\medskip

We may now prove the theorem. Clearly $K$ is non-commutative,
otherwise, $H_0\cap K=1$, because of the above proposition, and
we would then have $G=H_0.K$. This shows that $K$ contains at least one non-commutative
indecomposable direct factor.

Let $\mathcal{X}$ be the class up to isomorphism of an indecomposable non-commutative direct
factor of $K$ of minimal order. If $L$ is a direct factor of $G$
which is a member of $\mathcal{X}$ then $L'\subset
H_0'$. That is true because $L/(L\cap H_0')$ is isomorphic to a
direct factor of $G/H_0'\simeq (H_0/H_0')\times K$. But $L\cap H_0'>1$
according to the corollary of proposition 5.2, which shows that all
the indecomposable direct factors of $L/(L\cap H_0')$ have strictly less
elements than a member of $\mathcal{X}$.
By construction of $\mathcal{X}$, all
the indecomposable direct factors of $L/(L\cap H_0')$ are commutative,
so $L/(L\cap H_0')$ is itself commutative. We have shown that $L'\subseteq
H_0'$.

Let $N$ be a direct factor of $G$ isomorphic to a direct product of
members of $\mathcal{X}$. Also assume that $N$ is maximal in that respect.
It now suffices to prove that $N$ is isomorphic to a direct factor of $H_0$,
and we will have shown that a counter-example to our theorem cannot exist.

Clearly $N'\subseteq H_0'$. $H_0/H_0'$ is a subset of the center of
$G/H_0'$ because it is a commutative direct factor of that group. Likewise,
$(N.H_0')/H_0'$ is a commutative direct factor of $G/H_0'$ because
$N'\subseteq H_0'$. Therefore $(N.H_0)/H_0\subseteq
\mathcal{Z}(G/H_0')=(H_0/H_0').(((\mathcal{Z}(G).H_0)\cap
M)/H_0')$ (by lemma 4.2). So $H_0.N\subset
H_0.\mathcal{Z}(G)$.

Now $H_0.\mathcal{Z}(G)=H_0.S$, where $S$ is a supplementary of
$\mathcal{Z}(H_0)$ in $\mathcal{Z}(G)$. Clearly $H_0\cap S=1$ and
therefore $H_0$ is a direct factor of $H_0.\mathcal{Z}(G)\supset
H_0.N$. We have shown that $H_0$ is a direct factor of $H_0.N$
and that its supplementary is commutative. On the other hand, $N$ is a
direct factor of $G$ and thereby also of $H_0.N$. We now use the Remak-Krull-Schmidt
theorem on $H_0.N$. Notice that $N$ has
no non-trivial commutative direct factors to obtain that $N$ is
isomorphic to a direct factor of $H_0$. This is precisely what we have striven
to show. Our theorem is now proven.

\medskip

\section{Some additional remarks}
\begin{enumerate}
\item The theorem no longer holds if $G$ is infinite.
We give a simple counter-example. Let $G=(\Z/p.\Z)^{\N}\times (\Z/p^2.\Z)^{\N}$,
$H$ (resp. $K$) the subgroup of $G$ consisting of all pairs $(f,g)$ such that $f:\N
\longmapsto \Z/p.\Z$ and  $g:\N \longmapsto \Z/p^2.\Z$ with
$f(2n)=g(n)=0$ (resp. $f(2n+1)=0$).

Take $H_0$ the subgroup of $G$ consisting of all pairs $(f,g)$ with
$g(2n)=0$ and $g(2n+1)\in p.\Z/p^2.\Z$.

Clearly, $G=H.K$ and $H\cap K=1$. Moreover $H_0$ is not a direct factor
of $G$. But $H\simeq H_0\simeq (\Z/p.\Z)^{\N}$ and $K\simeq
G/H_0\simeq (\Z/p.\Z)^{\N}\times (\Z/p^2.\Z)^{\N}$.

\item The following statement does not hold: "Let $G$ be a finite group.
If there exists a split extension
$1\longmapsto H \longmapsto G \longmapsto K \longmapsto 1$ 
then any extension $1\longmapsto H \longmapsto G
\longmapsto K \longmapsto 1$ also splits".

We give a counter-example. Let $A$, $B$, $C$ and $D$ be four groups, each one
of them isomorphic to $\Z/p.\Z$, and $e_A$ a generator of $A$. $A$ acts
on $B\times C$ by $\phi: k.e_A\longmapsto
[(b,c)\longmapsto (b,c+k.b)]$.

Set $G=(A\ltimes_{\phi}(B\times C))\times D$. Clearly $B\times
C$ is a normal subgroup of $G$, with $A\times D$ a supplementary subgroup.
This defines a split extension $1\longmapsto
(\Z/p.\Z)^2\longmapsto G \longmapsto (\Z/p.\Z)^2\longmapsto 1$.

However $C\times D$ is in the center of $G$ and $G/(C\times D)$ is
isomorphic to $(\Z/p.\Z)^2$. On the other hand $C\times D$ cannot
have a supplementary because as it is in the center, the semi-direct product
would be trivial and $G$ and $(\Z/p.\Z)^4$ would be isomorphic.
So $1\longmapsto
(\Z/p.\Z)^2\longmapsto G \longmapsto (\Z/p.\Z)^2\longmapsto 1$ is an extension
which does not split.

\end{enumerate}

\medskip

\medskip

BIBLIOGRAPHY:

\begin{enumerate}
\item Group Theory I: \emph{Michio Suzuki}
\item Group Theory : \emph{W.R. Scott}
\item Maximal Subgroups of Direct Products: \emph{Jacques ThÈvenaz}
\end{enumerate}

\end{document}